%7.12.07
%4gvuj
\input amstex\input amsppt.sty
\magnification=\magstep1
%\hsize=30truecc\baselineskip=16truept
\voffset=-1cm
\advance\vsize-1cm\hoffset0cm\advance\hsize1cm
\NoBlackBoxes
\define\R{{\Bbb R}}\define\Z{{\Bbb Z}}
\def\pr{\mathop{\fam0 pr}}

\def\Emb{\mathop{\fam0 Emb}}

\def\id{\mathop{\fam0 id}}

\def\lk{\mathop{\fam0 lk}}
\def\delet{\mathaccent"7017 }
\def\rel#1{\allowbreak\mkern8mu{\fam0rel}\,\,#1}
\def\Int{\mathop{\fam0 Int}}

\def\Cl{\mathop{\fam0 Cl}}
\def\Emb{\mathop{\fam0 Emb}}

\def\link{\mathop{\fam0 link}}
\def\im{\mathop{\fam0 im}}

\def\t{\widetilde}

%In particular, please check the 2000 Mathematics Subject Classification
%- Page 10, line 4:  text slightly modified. Please check.
%change
%- Page 10, line -11: "define" - ? maybe you mean "decompose".
%no
%- Page 16:  [MR05] updated.
%ok, add arxiv

\topmatter
\title A new invariant and parametric connected sum of embeddings 
\endtitle
\author A. Skopenkov \endauthor
\address Department of Differential Geometry, Faculty of Mechanics and
Mathematics, Moscow State University, Moscow, Russia 119992, and Independent
University of Moscow, B. Vlas\-yev\-skiy, 11, 119002, Moscow, Russia.
e-mail: skopenko\@mccme.ru \endaddress
\thanks 
The author gratefully acknowledges the support by INTAS Grant No. YSF-2002-393, 
by the Russian Foundation for Basic Research, Grants No 05-01-00993, 
07-01-00648a and 06-01-72551-NCNILa, 
President of Russian Federation Grants MD-3938.2005.1, MD-4729.2007.1 and 
NSH-4578.2006.1, and by 
the Pierre Deligne fund based on his 2004 Balzan prize in mathematics.
\endthanks

\abstract
We define {\it an isotopy invariant} of embeddings $N\to\R^m$ of manifolds into 
Euclidean space. 
This invariant together with the $\alpha$-invariant of Haefliger-Wu is complete 
in the dimension range where the $\alpha$-invariant could be incomplete. 
We also define {\it parametric connected sum} of certain embeddings 
(analogous to surgery). 
This allows to obtain new completeness results for the $\alpha$-invariant
and the following estimation of isotopy classes of embeddings. 
{\it For the piecewise-linear category, a $(3n-2m+2)$-connected $n$-manifold 
$N$ and $\frac{4n+5}3\le m\le\frac{3n+2}2$ each preimage of $\alpha$-invariant 
injects into 
%are in 1--1 correspondences with 
a quotient of $H_{3n-2m+3}(N)$, where the coefficients are $\Z$ for $m-n$ odd 
and $\Z_2$ for $m-n$ even.} 
\endabstract
\endtopmatter

\document
\hfill{\it Dedicated to the centennary of K. Borsuk}

\head 1. Introduction and main results 
\footnote{This paper is to appear in Fundamenta Mathematicae 197 (2007).} 
\endhead 

This paper is on the classical Knotting Problem: {\it for an $n$-manifold $N$ 
and a number $m$ describe the set $\Emb^m(N)$ of isotopy classes of embeddings 
$N\to\R^m$}.
For recent surveys see [RS99, Sk07]; whenever possible we refer to these 
surveys not to original papers. 

%The Knotting Problem is more accessible for $2m\ge3n+4$.  
%%For $2m\ge3n+4$ the set $\Emb^m(N)$ was described using the {\it Haefliger-Wu
%%invariant} and the invariants that could be derived from it. For many cases 
%%this description allows explicit calculations [Ha63, RS99, \S4, Sk02, \S1, 
%Sk07, \S5].
%The Knotting Problem is much harder for $2m\le3n+3:$
%if $N$ is a closed manifold that is not a disjoint union of spheres, then no 
%explicit complete descriptions of isotopy classes was known

All known complete concrete classification results (except for the Haefliger
classification of links and smooth knots and recent results [KS05, Sk06, Sk06', 
CRS07, CRS]) can be obtained using $\alpha$-invariant of Haefliger-Wu (defined 
below).
For another approaches see [Br68, GW99, CRS04, We].
%in spite of the existence of interesting approaches [Br68, Wa70, GW99]

We define an isotopy invariant of embeddings which, together with the 
$\alpha$-invariant, is complete in the dimension range where the 
$\alpha$-invariant could be incomplete (the $\beta$-invariant Theorem of \S2).
We also define {\it parametric connected sum} of certain embeddings (see the 
end of \S2; this is a 'surgery' of an embedding preserving the embedded 
manifold). 
This allows to obtain new estimations of isotopy classes of embeddings and 
completeness results for the $\alpha$-invariant (the New Isotopy and Embedding 
Theorems of \S1).

We work in the piecewise linear (PL) category [RS72]. 
(By [Br72] for $m\ge n+3$ the classification of embeddings of PL manifolds is 
the same in the PL and the TOP categories. 
Analogously to [Sk06'] our results give some information for the smooth 
category.) 

%\bigskip
%{\bf Main results.} 
%The {\it Haefliger-Wu invariant} defined as follows.

Let
$$\t N=\{(x,y)\in N\times N\ |\ x\ne y\}$$
be the {\it deleted product} of $N$, i.e. the configuration space
of ordered pairs of distinct points of $N$.
For an embedding $f:N\to\R^m$ one can define a map $\t f:\t N\to S^{m-1}$
by the Gauss formula
$$\t f(x,y)=\frac{fx-fy}{|fx-fy|}.$$
This map is equivariant with respect to the 'exchange of factors' involution
$t(x,y)=(y,x)$ on~$\t N$ and the antipodal involution on~$S^{m-1}$.

Define $\alpha(f)$ to be the
equivariant homotopy class of the map $\t f$, cf. [Gr86, 2.1.E].
This is clearly an isotopy invariant.

%Let $\Emb^m(N)$ be the set of PL embeddings $N\to\R^m$ up to PL isotopy.
Let $\pi^{m-1}_{eq}(\t N)$ be the set of equivariant maps $\t N\to S^{m-1}$
up to equivariant homotopy.
Thus the $\alpha$-invariant is a map
$$\alpha:\Emb\phantom{}^m(N)\to\pi^{m-1}_{eq}(\t N).$$
It is important that using algebraic topology methods the set
$\pi^{m-1}_{eq}(\t N)$ can be explicitly calculated in many cases
[BG71, Ba75, Ya83, 
RS99, Sk02, GS06, Sk07, \S5]. 
So it is very interesting to know under which conditions the $\alpha$-invariant 
is bijective.

\smallskip
{\bf Isotopy Theorem.}
{\it (a) [Sk07, the Haefliger-Weber Theorem 5.4]
The $\alpha$-invariant is bijective for embeddings $N\to\R^m$ of an
$n$-polyhedron $N$, if
$$2m\ge3n+4.$$
\quad (b) [Sk07, Theorem 5.5]
The $\alpha$-invariant is bijective for $m\ge n+3$ and embeddings
$N\to\R^m$ of a closed $k$-connected $n$-manifold $N$, if  
$$2m\ge3n+3-k.$$} 

These theorems have many specific corollaries [Sk07].

In this paper we study the case one dimension lower than in the Isotopy Theorem 
(b).

By $\Z_{(d)}$ we denote $\Z$ for $d$ even and $\Z_2$ for $d$ odd. 
%If the coefficients of a homology group are omitted, then they are $\Z$. 

\smallskip
{\bf New Isotopy Theorem.}
{\it Let $N$ be a closed $k$-connected orientable $n$-manifold, 
$$2m=3n+2-k\quad\text{and}\quad n\ge3k+6\ge6.$$  
%\frac{4n+5}3\le m$.
\quad (a) The $\alpha$-invariant is surjective and each its preimage maps 
injectively into certain quotient of $H_{k+1}(N;\Z_{(m-n-1)})$. 

(b) The $\alpha$-invariant is bijective if either 
$(n,k,m)=(6,0,10)$ or $N$ is almost parallelizable and $(n,k,m)=(n,n-14,n+8)$, 
where $14\le n\le18$.} 

\smallskip
The new part of the New Isotopy Theorem is {\it estimation of point preimages} 
of the $\alpha$-invariant (which is surjective by the Embedding Theorem (b) 
below). 
These preimages could apriori be non-trivial by [Sk06', Example 1.6.b] stated 
below, and could depend on $n$, $k$, $N$ and the element of 
$\pi^{m-1}_{eq}(\t N)$ (of which we take the preimage).

For $m-n$ even the New Isotopy Theorem (a) implies that these preimages are
{\it finite}; the orientability assumption can be dropped.

The case $(n,k,m)=(6,0,10)$ of the New Isotopy Theorem (b) shows that 
[Ba75, Proposition 4] is true in the PL category for 6-manifolds. 

Under the assumptions of the New Isotopy Theorem the $\alpha$-invariant is not 
always injective:

{\it For each even $n\not\in\{6,14\}$ and $2m=3n+2$ the $\alpha$-invariant is 
not injective for embedings $S^1\times S^{n-1}\to\R^m$}
[Sk06', Example 1.6.b].
\footnote{Another examples of non-injectivity of the $\alpha$-invariant are
recalled in [RS99, \S4, Sk02, \S1, Sk07, \S5].} 

Some classification results for $(3n-2m+2)$-connected manifold 
$N=S^p\times S^q$ are obtained in [Sk06', Theorems 1.3 and 1.4, CRS07, CRS].
It is very surprising that something can be proved for {\it general} manifolds 
$N$.

\smallskip
{\bf Conjecture.} {\it (a) If $n\ge3k+4$ and $N$ is a closed $k$-connected 
almost parallelizable $n$-manifold, then there is an exact sequence of sets 
with an action $w$  
$$H_{k+1}(N,\Z_{(m-n-1)})\overset w\to\to\Emb\phantom{}^m(N)\overset
\alpha\to\to\pi^{m-1}_{eq}(\t N)\to0.$$ 
\quad (b) The Isotopy Theorem (a) holds for $(n,k,m)=(7,1,11)$ and $N$ spin, 
as well as for $(n,k,m)=(19,5,27)$ and $N$ almost parallelizable. 
\footnote{This follows from our proof of the Isotopy Theorem (a) (\S2) and an 
improvement [Sk06', Standardization Lemma] of the Standardization Lemma of \S2.} 
}

%BG71 nichego ne daet, t.k. nuzhno n\le3k

\smallskip
The corresponding known and new {\it surjectivity} results are as follows.
\footnote{Examples of non-surjectivity of the $\alpha$-invariant are 
recalled in [RS99, \S4, Sk07, \S5].} 

\smallskip
{\bf Embedding Theorem.}
{\it (a) The $\alpha$-invariant is surjective for embeddings $N\to\R^m$ of 
an $n$-polyhedron $N$, if $2m\ge3n+3$ [Sk07, the Haefliger-Weber Theorem 5.4].

(b) The $\alpha$-invariant is surjective for embeddings $N\to\R^m$ of a
closed $k$-connected $n$-manifold $N$ when $2m\ge3n+2-k$ and $m\ge n+3$ 
[Sk07, Theorem 5.5].}

\smallskip
{\bf New Embedding Theorem.}
{\it Let $N$ be a closed $k$-connected $n$-manifold. 
% and $m=(3n+1-k)/2$. Assume that $n\ge3k+5\ge5$ and $m-n$ is 4 or 8.
%$\frac{4n+4}3\le m\le\frac{3n+1}2$ 
The manifold $N$ embeds into $\R^m$ if there is an equivariant map 
$\t N\to S^{m-1}$ and either  

$(n,k,m)=(7,0,11)$ and $N$ is orientable, or 

$(n,k,m)=(8,1,12)$ and $N$ is spin, or 

$(n,k,m)=(n,n-15,n+8)$ and $N$ is almost parallelizable, where $15\le n\le20$.}

%The pairs $(n,m,d)$ satisfying to the assumptions of the New 
%Embedding Theorem are
%$(7,11,0)$, $(8,12,1)$ and $(n,n+8,n-15)$, where $15\le n\le20$.

\smallskip
An $n$-manifold is {\it $p$-parallelizable} if any embedding $S^p\to N$ 
can be extended to an embedding $S^p\times D^{n-p}\to N$.
Note that 1-parallelizability is equivalent to orientability and 
1\&2-parallelizability is equivalent to being a spin manifold. 
\footnote{It would be interesting to reformulate the 
$(k+1)$-parallelizability condition 
for $k$-connected manifolds in terms of Stiefel-Whitney classes.} 
The almost parallelizability condition in the results of \S1 can be relaxed to 
$(k+1)$-parallelizability. 

\bigskip
{\bf Acknowledgements.} 

These results were presented at the Borsuk Centennary Conference (Bedlewo, 
2005). 
%and announced in [Sk05]. 
I would like to acknowledge M. Skopenkov and S. Melikhov for useful discussions.

%\bigskip
\newpage
\head 2. Proofs \endhead 

%\bigskip
%\smallskip
{\bf Almost embeddings and almost concordances.}

%Recall that we work in the PL category. 

An embedding $F:N\times I\to\R^m\times I$ is a {\it concordance} if
$N\times0=F^{-1}(\R^m\times0)$ and $N\times1=F^{-1}(\R^m\times1)$.
We tacitly use the facts that in codimension at least 3

{\it concordance implies isotopy} [Hu70, Li65], and

{\it every concordance or isotopy is ambient} [Hu66, Ak69].%, Edw75].

Let $N$ be a connected $n$-manifold and $B^n\subset\delet N$ some $n$-ball.
The {\it self-intersection set} of a map $F:N\to\R^m$ is
$$\Sigma(F):=\Cl\{x\in N\ |\ \#F^{-1}Fx\ge1\}.$$
A map $F:N\to\R^m$ is an {\it almost embedding} of $(N,B^n)$, if
$\Sigma(F)\subset B^n$.
Or, equivalently, if $F|_{N-B^n}$ is an embedding and
$F(N-B^n)\cap F(B^n)=\emptyset$.

A map $F:N\times I\to\R^m\times I$ is an {\it almost concordance}
of $(N,B^n)$ if
$$N\times0=F^{-1}(\R^m\times0),\quad N\times1=F^{-1}(\R^m\times1)\quad
\text{and}\quad\Sigma(F)\subset B^n\times I.$$
Instead of the pair $(N,B^n)$ we shall always write simply $N$
(no confusion would arise).
\footnote{Almost embeddings and almost concordances were called 
quasi-embeddings and quasi-con\-cor\-dan\-ces in [Sk02].}
\footnote{
Fix points $x=+1\in S^0\subset S^p$ and $y=+1\in S^0\subset S^{n-p}$.
By general position for $m\ge n+p+1$ any map $f:S^p\times S^{n-p}\to\R^m$ such 
that $\Sigma(f)\cap x\times S^{n-p}=\emptyset$ is homotopic through such maps 
to a map $f'$ whose self-intersection set is contained in the ball 
$D^p_-\times D^{n-p}_-$. 
Analogous statement holds for $m\ge n+p+2$ and almost concordances.
In this sense for $m\ge n+p+2$ 
%and $p\ge1$??? 
the above definition of an almost concordance for $N=S^p\times S^{n-p}$ and 
$B^n=D^p_-\times D^{n-p}_-$ agrees with that of [Sk06', \S2].}

\smallskip
{\bf Almost Embeddings Theorem.}
{\it Suppose that $N$ is a closed $k$-connected $n$-manifold, $k\ge0$ and 
$m\ge n+2$. 

(a) If $f,g:N\to\R^m$ are almost concordant embeddings,  
then $\alpha(f)=\alpha(g)$ [Sk02, Theorem 5.2.$\alpha$].

% (and, moreover, $\alpha_G(f)=\alpha_G(g)$ for each $G$

(b) If $2m=3n+2-k$ and $f,g:N\to\R^m$ are embeddings such that
$\alpha(f)=\alpha(g)$, then $f$ and $g$ are almost concordant 
[Sk02, Theorem 2.2.q]. 

(c) If $2m=3n+1-k$ and $\varphi\in\pi^{m-1}_{eq}(\t N)$, then there is 
an almost embedding $F:N\to\R^m$ such that $\alpha(F)=\varphi$ 
[Sk02, Theorem 2.2.q].}

\bigskip
{\bf Appendix: some results and conjectures on almost embeddings.}

This section is not used in the proof of main results, but is perhaps of 
independent interest. 

Complete classification of embeddings of a given $n$-manifold $N$ into
$S^{n+2}$ up to isotopy (or concordance) seems to be hopeless because it is
such for $N=S^n$.
So it is interesting to obtain complete classification of embeddings of a
given $n$-manifold $N$ into $S^{n+2}$ 'modulo knots $S^n\to S^{n+2}$'.
The notion of almost concordance is not only useful to study the initial
problem (of classification of embeddings up to concordance) for $m\ge n+3$, but
is a good notion of 'concordance modulo knots $S^n\to S^{n+2}$', because any
knot $S^n\to S^{n+2}$ is almost concordant to the trivial knot.
\footnote{For $N=S^{n_1}\sqcup\dots\sqcup S^{n_k}$ the classification of 
embeddings $N\to\R^{n+2}$ up to {\it link homotopy} is motivated by 
classification of embeddings $N\to\R^{n+2}$ 'modulo knots $S^n\to\R^{n+2}$'.}
Cf. [MR05]. 

We conjecture that almost concordance is equivalent to another natural
equivalence relation of 'concordance modulo knots $S^n\to S^{n+2}$', i.e. that
for a closed $n$-manifold $N$ two embeddings $N\to S^{n+2}$
are almost concordant if and only if one can be obtained from an embedding
concordant to the other by connected summation with knots $S^n\to S^{n+2}$.

Parts (a) and (b) of the following corollary are implied by the Almost 
Embeddings Theorem (b) and by [Sk02, Theorem 2.3.q], respectively.

\smallskip
{\bf Corollary.}
{\it (a) Let $N$ be a sphere with $g$ handles. 
Then the set of PL almost embeddings $N\to\R^4$ up to PL almost concordance is 
in 1--1 correspondence with $\Z_2^{2g}\cong H_1(N;\Z_2)$.

(b) For a closed simply-connected 4-manifold $N$, the set of smooth
almost embeddings $N\to\R^6$ up to PL almost concordance is in 1--1
correspondence with $\pi^5_{eq}(\t N)$.}

\smallskip
We conjecture that the set of PL embeddings $S^1\times S^1\to\R^4$ up to PL
almost concordance consists of exactly 3 elements (i.e. that the almost
embedding $S^1\times S^1\to\R^4$, corresponding by Corollary (a) to the class
$(1,1)\in H_1(S^1\times S^1;\Z_2)$, is not almost concordant to a PL
embedding).
\footnote{A related result states that any PL embeding  
$S^2\sqcup\dots\sqcup S^2\to S^4$
is link homotopic to the trivial embedding [BT99, Ba01] (the case of two
components was proved earlier by Hosokawa-Suzuki).} 

The restriction $k\ge0$ is essential in the Almost Embeddings Theorem (b).
\footnote{Indeed, for $l\not\in\{3,7\}$ take a link 
$f:S^0\times S^{2l-1}\to\R^{3l}$
such that $\lambda_{12}(f)=\lambda_{21}(f)=[\iota_l,\iota_l]$.
Then $\alpha(f)=\Sigma^\infty\lambda_{21}(f)=0$ but $f$ is not almost
concordant to the standard embedding.} 

We conjecture that the Almost Embeddings Theorem (b) holds in the smooth 
category, and that in [Sk02, Theorem 2.3.q] and in the injectivity of
[Sk02, Theorem 2.3.$\alpha$] we can replace PL category by DIFF (if $N$ is 
a smooth manifold).
\footnote{This could perhaps be proved analogously to the cited results 
using the relative version of [Sk02, Disjunction Theorem 3.1].} 

%For a manifold $N$ with a fixed decomposition $N=N_1\sqcup\dots\sqcup N_k$
%a map $f:N\to\R^m$ is called a {\it almost-embedding}, if
%$fN_i\cap fN_j=\emptyset$ for $i\ne j$ and $\Sigma(f|_{N_i})$ is contained
%in some ball in $\delet N_i$.
%This concept generalizes both link maps and almost-embeddings of connected
%manifolds.

We conjecture that for $n$ even, $m=(3n+1-k)/2\ge n+3$ and a $k$-connected 
closed $n$-manifold $N$ such that $H_{k+1}(N)$ is free there is an exact 
sequence of sets
$\Emb\phantom{}^m(N)\overset\alpha\to\to\pi\phantom{}^{m-1}_{eq}(\t N)
\overset\beta\to\to H_{k+1}(N)$. 
Cf. the New Embedding Theorem.
\footnote{By [Sk02, Theorem 2.3.q] and the $\beta$-invariant Theorem
it suffices to prove that for an almost embedding $F:N\to\R^m$ the obstruction
$\beta(F)$ does not depend on $F|_{B^n}$ (only on $F|_{N-\delet B^n}$).
This obstruction is a map $b:\pi_n(M)\to H_{k+1}(N)$.
We can prove that this map is constant by checking that $\pi_n(M)$ is finite 
and for fixed $\varphi_0\in\pi_n(M)$ the map 
$\varphi\mapsto b(\varphi)-b(\varphi_0)$ is a homomorphism.}

\bigskip
{\bf Definition of $\beta$-invariant.}

The Almost Embeddings Theorem (b) suggests the definition of an invariant,
required for classification of embeddings when $2m\le3n+2-k$.
For each almost concordance $F$ between embeddings analogously to [Hu69, 
XI.4.iii, Hu70', p. 408, Ha84, \S1] we define an obstruction $\beta(F)$ to 
modification of $F$ to a concordance.
Roughly speaking, $\beta(F)$ measures the linking of $\Sigma(F)$ with $F(N)$.

Analogous invariants are the Sato-Levine invariant of knots, the 
Hudson-Habegger obstruction to embedding disks and Fenn-Rolfsen-Koschorke-Kirk 
$\beta$-invariant of link maps, see references in [Sk06'].
In the proof of the New Embedding and Isotopy Theorems we do not use the
definition but only use the properties of $\beta$-invariant
(they are stated in the $\beta$-invariant Theorem of the next subsection).

In this and the next subsections we omit the coefficients $\Z_{(m-n-1)}$ of 
chain groups in the notation. 

Suppose that $N$ is a connected orientable $n$-manifold (possibly, with
boundary) and $F:N\to B^m$ is a proper general position almost embedding whose
restriction to the boundary is an embedding. 
($F$ could be an almost concordance 
between embeddings.) 
Take a triangulation $T$ of $N$ such that $B^n$ is a subcomplex of
$T$ and $F$ is linear on simplices of $T$.
Then $\Sigma(F)$ is a subcomplex of $T$.
Denote by $[\Sigma(F)]\in C_{2n-m}(B^n)$ the sum of
top-dimensional simplices of $\Sigma(F)$.

For $m-n$ odd the coefficient $\pm1$ of an oriented simplex
$\sigma\subset\Sigma(F)$ is defined as follows. 
\footnote{For $m-n$ even this sign can also be defined but is not used.}
Fix in advance any orientation of $N$ and of $B^m$.
By general position there is a unique
simplex $\sigma'$ of $T$ such that $F(\sigma)=F(\sigma')$.
The orientation on $\sigma$ induces an orientation on $F\sigma$ and then on
$\sigma'$.
The orientations on $\sigma$ and $\sigma'$ induce orientations
on normal spaces in $N$ to these simplices.
These two orientations (in this order) together with the orientation on
$F\sigma$ induce an orientation on $B^m$.
If this orientation agrees with the fixed orientation of $B^m$,
then the coefficient of $\sigma$ is $+1$, otherwise $-1$.
Clearly, change of orientation of $\sigma$ changes the sign of $\sigma$ in
$[\Sigma(F)]$, so the sign is well-defined.
\footnote{This definition of sign is equivalent to Hudson's one given as follows.
The orientation on $\sigma$ induces an orientation on $F\sigma$ and on $\sigma'$,
hence it induces orientation on their links.
Consider the oriented $(2m-2n-1)$-sphere $\lk F\sigma$, that is the link
of $F\sigma$ in certain triangulation of $B^m$, 'compatible' with $T$.
This sphere contains disjoint oriented $(m-n-1)$-spheres $F(\lk_T\sigma)$
and $F(\lk_T\sigma')$.
Their linking coefficient
$\link_{\lk F\sigma}(F(\lk_T\sigma),F(\lk_T\sigma'))\in\Z_{(m-n-1)}$ is the
coefficient of $\sigma$ in $[\Sigma(F)]$, which equals $\pm1$.}

By [Hu69, Lemma 11.4, Hu70', Lemma 1] $\partial[\Sigma(F)]=0$. 
\footnote{Here we use the fact that coefficients are $\Z_2$ for $m-n$ even.}
Hence 
$$[\Sigma(F)]=\partial C\quad\text{for some}\quad C\in C_{2n-m+1}(B^n).$$
By [Hu70', Corollary 1.1] $\partial FC=0$.
Hence 
$$F(C)=\partial D\quad\text{for some}\quad D\in C_{2n-m+2}(B^m).$$
By general position 
$$\widetilde D:=(F|_{N-\Int B^n})^{-1}(D)\in C_{3n-2m+2}(N).$$ 
Since the support of $C$ is in $B^n$, the support of $F(C)=\partial D$ is in 
$F(B^n)$.
Hence $\widetilde D$ is a cycle and we can define 
\footnote{This definition agrees with that for $N=S^p\times S^q$ [Sk06', subsection
'a new embedding invariant' of \S2] when $m\ge 2p+q+2$.} 
$$\beta(F):=[\widetilde D]\in H_{3n-2m+2}(N;\Z_{(m-n-1)}).$$

%Denote $N_0:=N-\Int B^n$

%\smallskip
{\it Proof that $\beta(F)$ is well-defined, i.e. is independent of the choices
of $C$ and $D$.} 
%\footnote{Analogous to the particular case $N=S^p\times S^q$ [Sk06', \S2].}
The independence of the choice of $D$ is standard. 
Let us prove the independence of the choice of $C$.  
For an almost embedding $F$ if $\partial C_1=\partial C_2=[\Sigma(F)]$ 
then $\partial(C_1-C_2)=0$. 
Hence 
$$C_1-C_2=\partial X\quad\text{for some}\quad X\in C_{2n-m+1}(B^n).
$$ 
Thus $FC_1-FC_2=\partial FX$. 
Hence we can take chains $D_1$ and $D_2$ as above and such that $D_1-D_2=FX$. 
Since the support of $X$ is in $B^n$, we have $\widetilde D_1=\widetilde D_2$. 
\qed

\bigskip
{\bf Properties of the $\beta$-invariant and proof of the New Isotopy Theorem 
(a).}

\smallskip
{\bf $\beta$-invariant Theorem.}
{\it Let $N$ be a connected orientable $n$-manifold (possibly with boundary) 
and $m\ge n+3$.
To each proper general position almost embedding $F:N\to B^m$ whose restriction 
to the boundary is an embedding there corresponds an element
$$\beta(F)\in H_{3n-2m+3}(N;\Z_{(m-n-1)})\qquad\text{such that}$$

(obstruction) If $F$ is an embedding, then $\beta(F)=0$.

(invariance) $\beta(F)$ is invariant under almost concordance of
$F$ relative to the boundary.

(completeness) If $\beta(F)=0$ and $N$ is homologically $(3n-2m+1)$-connected,
then $F$ is almost concordant $\rel(N-\delet B^n)$ to an embedding.

(additivity)  Suppose that $N=X\times I$ and $F$, $F'$ are almost concordances  
between $f_0$ and $f_1$, $f_1$ and $f_2$, respectively, 
Denote by $\overline F$ the reversed $F$, i.e. $\overline F(x,t)=F(x,1-t)$.
Define an almost concordance $F\cup F'$ between $f_0$ and $f_2$ as 'the union' 
\footnote{Observe that the union of almost concordances is associative up to ambient
isotopy.}
of
$$F:X\times[0,1]\to\R^m\times[0,1]\quad\text{and}
\quad F':X\times[1,2]\to\R^m\times[1,2].$$
Then
$\beta(F\cup F')=\beta(F)+\beta(F')$ and $\beta(\overline F)=-\beta(F)$.}

\smallskip
Here the orientability assumption can be dropped for $m-n$ even.

The obstruction and additivity follow obviously by
the definition of $\beta$-invariant. 

The invariance is analogous to [Hu70', Lemma 2, cf. Hu69, Lemma 11.6]. 
The completeness is a non-trivial property, but it is an easy consequence 
of known results [Hu70', Theorem 2, Ha84, Theorem 4]. 
See the details below.

\smallskip
{\it Proof of the invariance.}
Let $F_0$ and $F_1$ be two almost concordances between embeddings $f$ and $g$.
Suppose that $\Phi:N\times I\to B^m\times I$ is an almost concordance between
$F_0$ and $F_1$.
%Then $\Sigma(\Phi)\subset B^n\times I$.
As in the above definition of $\beta(F_0)$ and $\beta(F_1)$ take chains
$$C_0,C_1\in C_{2n-m+1}(B^n)\quad\text{such that}\quad
\partial C_j=[\Sigma(F_j)]\quad\text{for}\quad j\in\{0,1\}.$$
Analogously to the above definition of $\beta(F)$ define chain
$$[\Sigma(\Phi)]\in C_{2n-m+1}(B^n\times I).$$
Let $i_j:N\cong N\times j\to N\times I$ be the inclusions.
Then
$$\partial[\Sigma(\Phi)]=i_1[\Sigma(F_1)]-i_0[\Sigma(F_0)],\quad\text{so}\quad
\partial([\Sigma(\Phi)]+C_0-C_1)=0.$$
Therefore there exists
$$C\in C_{2n-m+2}(B^n\times I)\quad\text{such that}
\quad\partial C=[\Sigma(\Phi)]+C_1-C_0.$$
Analogously to [Hu70', Lemma 2] $\partial\Phi C=i_1(F_1C_1)-i_0(F_0C_0)$.
Let $\pr:N\times I\to N$ be the projection.
Hence $\partial\pr\Phi C=F_1C_1-F_0C_0$. 
Thus analogously to the proof of the independence of $\beta$ of the choice of 
$C$ we obtain $\beta(F_0)=\beta(F_1)$. 
\qed

\smallskip
{\it Proof of the completeness.}
Denote
$$M:=B^m-\Int R_{B^m}(F(N-\delet B^n),F\partial B^n).$$
Observe that $F|_{B^n}:B^n\to M$ is a proper map whose restriction to the
boundary is an embedding.

Consider the following composition of Alexander and Poincar\'e duality 
isomorphisms (with the $\Z$-coefficients)
$$H_i(M)\cong H^{m-i-1}(B^m-M,\partial B^m-M)
\cong H^{m-i-1}(N-\delet B^n,\partial N)\cong$$
$$\cong H_{i+n-m+1}(N-\delet B^n,\partial B^n)\cong H_{i+n-m+1}(N).$$
Since $N$ is homologically $(3n-2m+1)$-connected, we obtain that $M$ is
homologically $(2n-m)$-connected.
Since $M$ is simply-connected, it follows that $M$ is $(2n-m)$-connected.
Then by [Hu70', Theorem 2, Ha84, Theorem 4] the class 
$[FC]\in H_{2n-m+1}(M;\Z_{(m-n-1)})$ is the complete obstruction to
the existence of a homotopy $\rel\partial B^n$ from $F|_{B^n}:B^n\to M$ to an
embedding.
This class goes to $\beta(F)=0$ under the composition of the above isomorphisms 
with coefficients $\Z_{(m-n-1)}$. 
Hence $F|_{B^n}$ is homotopic $\rel\partial B^n$ to an embedding.
Extending this embedding over $N$ by $F$ we obtain the required embedding 
$N\to B^m$.
\qed

\smallskip
{\it Proof of the New Isotopy Theorem (a).}
Fix any $\varphi\in\pi^{m-1}_{eq}(\t N)$ and any embedding $f_0:N\to\R^m$ such
that $\alpha(f_0)=\varphi$.
Define
$$K:=\{\beta(F_0)\in H_{k+1}(N;\Z_{(m-n-1)})\ |
\ \text{ $F_0$ is an almost concordance from $f_0$ to $f_0$}\}.$$
By the additivity of $\beta$-invariant, $K$ is a subgroup (depending on 
$n,k,N,\varphi$).

For any embedding $f:N\to\R^m$ such that $\alpha(f)=\alpha(f_0)$ by Almost
Embeddings Theorem (b) there is an almost concordance $F$ from $f$ to $f_0$.
(This together with the additivity of $\beta$-invariant implies that $K$ does 
not depend on the choice of $f_0$.)
So we can define a map
$$B:\alpha^{-1}(\varphi)\to H_{k+1}(N;\Z_{(m-n-1)})/K
\quad\text{by}\quad B(f):=\beta(F)+K.$$
If $F$ and $F'$ are two almost concordances from  $f$ to $f_0$, then
$F'\cup\overline F$ is an almost concordance from $f_0$ to $f_0$.
Hence the map $B$ is well-defined by the additivity of $\beta$-invariant.

If $B(f)\in K$, then $\beta(F)=\beta(F_0)$ for some almost
concordance $F_0$ from $f_0$ to $f_0$. Then $F\cup\overline F_0$
is an almost concordance from $f$ to $f_0$, and by the additivity
of $\beta$-invariant $\beta(F\cup\overline F_0)=0$. Hence by the
completeness of $\beta$-invariant $f$ is concordant to $f_0$. 
Thus $B$ is injective. \qed

\bigskip
{\bf Parametric connected sum of embeddings.} 

By $S^p=D^p_+\bigcup\limits_{\partial D^p_+=S^{p-1}=\partial D^p_-}D^p_-$ we 
denote the standard decomposition of $S^p$.
Analogously define $\R^m_\pm$ and $\R^{m-1}$.
Identify $D^p$ with $D^p_+$. 

For $m\ge n+2$ denote by $i$ the {\it standard embedding} which is the 
composition 
$S^p\times S^{n-p}\to\R^{p+1}\times\R^{n-p+1}\subset\R^m\subset S^m$.

Let $N$ be a closed connected $n$-manifold. 
Denote by $s:S^p\times D^{n-p}\to N$ an embedding. 
For the ball $B^n\subset N$ from the definition of an almost embedding 
(concordance) assume that $\im s\cap B^n=\emptyset$. 

A map $f:N\to S^m$ is called {\it $s$-standardized} if 

$f\circ s:S^p\times D^{n-p}\to D^m_-$ is the 
restriction of the standard embedding and 

$f(N-\im s)\subset\Int D^m_+$.

Roughly speaking, a map $N\to D^m$ is $s$-standardized if its image is put on 
hyperplane $D^{m-1}$ so that the image intersects the hyperplane in a 
standardly embedded $S^p\times D^{n-p}$ (indeed, for such a map the set $\im s$ 
can be pulled below the hyperplane to obtain an $s$-standardized embedding in 
the above sense). 
\footnote{Note that standardized in the sense of [Sk06', \S2] is 
$i$-standardized in the sense of this paper.} 

A concordance $F:N\times I\to S^m\times I$ between $s$-standardized maps
is called {\it $s$-standardized} if 

$F(\im s\times I)\subset S^m_-\times I$ is
the identical concordance and 

$F((N-\im s)\times I)\subset\Int S^m_+\times I$.

\smallskip
{\bf Standardization Lemma.}
{\it (a) If $m\ge n+p+2$, then any (almost) embedding $g:N\to S^m$ is isotopic 
to an $s$-standardized (almost) embedding.  

(b) If $m\ge n+p+3$, then any (almost) concordance between $s$-standardized 
embeddings $N\to S^m$ is isotopic relatively to the ends to an $s$-standardized 
(almost) concordance.}

%The proof is analogous to [Sk06', Standardization Lemma].

\smallskip
{\it Proof of (a).}
Fix a point $y\in D^{n-p}_-\subset S^{n-p}$.
Since $m\ge n+p+2\ge 2p+2$, it follows that $g|_{S^p\times y}$ is unknotted in
$S^m$.
So there is an embedding $\widehat g:D^{p+1}\to S^m$ such that
$$(*)\qquad\widehat g|_{\partial D^{p+1}}=g|_{S^p\times y}\quad\text{and}
\quad\widehat g\Int D^{p+1}\cap gN=\emptyset.$$ 
(The second property holds by general position because $m\ge n+p+2$.)
The regular neighborhood in $S^m$ of $\widehat gD^{p+1}$ is homeomorphic to 
the $m$-ball.
Take an isotopy moving this ball to $D^m_-$ and let $f':N\to S^m$ be the
embedding obtained from $g$.

Now we are done since the embedding $f'\circ s$ is isotopic to 
the standard embedding by the following result 
(because $m\ge n+3$, the pair $(S^p\times D^{n-p},S^p\times S^{n-p-1})$
is $(n-p-1)$-connected and $n-p-1\ge 2n-m+1$).
\qed

\smallskip
{\bf Unknotting Theorem Moving the Boundary.}
{\it Let $N$ be a compact $n$-dimensional PL manifold and $f,g:N\to D^m$
proper PL embeddings.
If $m\ge n+3$ and $(N,\partial N)$ is $(2n-m+1)$-connected, then $f$ and $g$
are properly isotopic} [Hu69, Theorem 10.2, p. 199].

\smallskip
{\it Proof of (b).}
This is a relative version of the proof of (a).
Take a concordance $G$ between standardized embeddings $f_0,f_1:N\to S^m$.
There is a level-preserving  embedding
$\widehat G:D^{p+1}\times\{0,1\}\to S^m\times\{0,1\}$ whose components satisfy 
to (*).
Since $m+1\ge p+1+3$, by the Haefliger-Zeeman Unknotting Theorem any 
concordance $S^p\times I\to S^m\times I$ standard on the boundary is isotopic 
to the standard concordance. 
Hence the map $\widehat G$ can be extended to an embedding
$\widehat G:D^{p+1}\times I\subset S^m\times I$
%\quad\text{
such that
%}\quad
$$\widehat G|_{\partial D^{p+1}\times I}=G|_{S^p\times y\times I}
\quad\text{and}\quad
\widehat G(\Int D^{p+1}\times I)\cap G(N\times I)=\emptyset.$$
(The second property holds by general position because $m\ge n+p+3$.)
Take a regular neighborhood
$$B^m\times I\quad\text{in}\quad S^m\times I\quad\text{of}\quad
\widehat GD^{p+1}\quad\text{such that}\quad
(B^m\times I)\cap(S^m\times\{0,1\})=D^m_-\times\{0,1\}.$$
Take an isotopy of $S^m\times I\rel S^m\times\{0,1\}$ moving $B^m\times I$ to
$D^m_-\times I$.
Let $F'$ be the concordance obtained from $G$ by this isotopy.

The embedding 
$F'|_{\im s\times I}:\im s\times I\to D^m_-\times I$ is 
isotopic $\rel D^m_-\times\{0,1\}$ to the identical concordance by the
following Unknotting Theorem Moving Part of the Boundary (which is proved
analogously to [Hu69, Theorem 10.2 on p. 199]).

{\it Let $N$ be a compact $n$-dimensional PL manifold, $A$ a codimension
zero submanifold of $\partial N$ and $f,g:N\to D^m$ proper PL embeddings.
If $m\ge n+3$ and $(N,A)$ is $(2n-m+1)$-connected, then $f$ and $g$
are properly isotopic $\rel\partial N-A$.}
\qed

%Let $T^n_p:=S^p\times S^{n-p}$ and 
%The sign $\sim$ between embeddings means that they are isotopic. 

\smallskip
By $R_k$ we denote the symmetry of $\R^k$ with respect to the plane 
$x_1=x_2=0$.

\smallskip
{\bf Definition of parametric connected sum.}
{\it (a) Take (almost) embeddings
$$f:N\to\R^m\quad\text{and}\quad g:S^p\times S^{n-p}\to\R^m.$$ 
If $m\ge n+p+2$, then by Standardization Lemma (a) we can make isotopies and 
assume that $f$ and $g$ are $s$-standardized and $i$-standardized, respectively.
Define an (almost) embedding 
$$f\#_sg:N\to\R^m\quad\text{by}\quad
(f\#_sg)(a)=\cases f(a) &a\not\in\im s\\
R_mg(x,R_{n-p}y)           &a=s(x,y)\endcases,$$ 
\quad (b) Take (almost) concordances 
$$F:N\times I\to\R^m\times I\quad\text{and}
\quad G:S^p\times S^{n-p}\times I\to\R^m\times I.$$ 
If $m\ge n+p+3$, then by the Standardization Lemma (b) we can make isotopies 
relative to the ends and assume that $F$ and $G$ are $s$-standardized and 
$i$-standardized, respectively.
Define a (almost) concordance
$F\#_sG:N\times I\to\R^m\times I$ by 
$$(F\#_sG)(a,t)=\cases F(a,t)& a\not\in\im s\\
(R_mG(x,R_{n-p}y,t),t)      &a=s(x,y)\endcases.$$} 

We do not need 
parametric connected sum to be independent on the choice of an almost 
concordance to a {\it standardized} almost embedding or almost concordance: we  
denote by $f\#_sg$ or $F\#_sG$ the result for {\it any} such choice.

%\smallskip
\bigskip
{\bf Proof of the New Isotopy Theorem (b) and the New Embedding Theorem.} 

\smallskip
{\bf The Hopf Invariant Lemma.} 
{\it Take the standard embedding $i:D^{p+1}\times S^q\to S^m$. 
Represent $\varphi\in\pi_{p+q}(S^{m-q-1})$ by a map (not necessarily an 
embedding) 
\linebreak
$\overline\varphi:S^{p+q}\to S^m-i(D^{p+1}\times S^q)\simeq S^{m-q-1}.$ 
If $2m=3q+2p+2$, then $\beta(i\#\overline\varphi)=\pm H\Sigma\varphi$. 
\footnote{Here $i\#\overline\varphi$ is embedded connected sum of linked 
embeddings but not parametrized connected sum as above; 
$i\#\overline\varphi=\overline\mu(\varphi)$ in the notation of [Sk06'].}
} 

\smallskip
{\it Proof (analogous to [Ko88, Theorem 4.8]).}   
We may assume that $\overline\varphi$ is a {\it smooth} general position framed 
immersion. 
%So in this proof we work in the smooth category
Extend $\overline\varphi$ to a {\it smooth} general position framed immersion 
$\widehat\varphi:B\to\R^m$, where $B:=B^{p+q+1}$. 
Then by [Ko88, Theorem 1.3] and [Ke59, Lemma 5.1] the class 
$\pm\Sigma\varphi\in\pi^S_{p+2q+1-m}$ is represented by the framed 
$(p+2q+1-m)$-submanifold 
$$\Delta:=\{(u,z)\in B\times S^q\ |\ \widehat\varphi(u)=i(a,z)\}\quad
\text{of}\quad B\times S^q\subset S^{p+2q+1-m}$$ 
(with natural framing). 
For a 0-chain $X$ with coefficients in $\Z_{(m-p-q-1)}$ in a connected manifold 
denote by $[X]$ the number of points in $X$ modulo 2 when $m-p-q$ is even (we 
need only this case for the Non-triviality Lemma) and the algebraic number of 
points when $m-p-q$ is odd. 
(I.e. $[X]$ is the 0-dimensional homology class of $X$.)  
Then by [Ko88, p. 411]    
$$\pm H\Sigma\varphi=
[\{(x,y)\in\Delta\times\Delta\ |\ \pr\phantom{}_2x=\pr\phantom{}_2y\}]$$ 
(this set is finite by general position). 
Thus  
$$\pm H\Sigma\varphi=[\{(u,v,z)\in B\times B\times S^q\ |
\ \widehat\varphi(u)=\widehat\varphi(v)=i(a,z)\}]=
[i(a\times S^q)\cap \widehat\varphi\pr\phantom{}_2 D]=
\beta(i\#\overline\varphi),$$ 
$$\text{where}\quad D:=\{(u,v)\in B\times B\ |
\ \widehat\varphi(u)=\widehat\varphi(v)\}.$$
Here the last equality holds because 
$$C:=\{(u,v)\in \partial B\times B\ |\ \overline\varphi(u)=\widehat\varphi(v)\},
\quad D\quad\text{and}\quad\widetilde D:=\widehat\varphi\pr\phantom{}_2 D$$ 
are as in the definition of $\beta$-invariant 
($C$ has natural orientation for $m-p-q$ odd); the groups 
$H_p(S^p\times S^q;\Z_{(m-p-q-1)})$ and $\Z_{(m-p-q-1)}$ are identified by the 
isomorphism $\gamma\mapsto\gamma\cap[a\times S^q]$. 
The definition of $\beta$-invariant does make 
sense in the piecewise-smooth category (and hence in the smooth category). 
Recall that piecewise-smooth category is equivalent to the PL category, i.e. 
the forgetful map from the set of PL embeddings up to PL isotopy to the 
set of piecewise differentiable embeddings up to piecewise differentiable 
isotopy is a 1--1 correspondence [Ha67].
\qed

\smallskip
{\bf Non-triviality Lemma.} {\it For $1\le p<l\in\{3,7\}$ there
exists an almost embedding $G:S^p\times S^{2l}\to\R^{3l+p+1}$ such
that $\beta(G)=1$.} 

\smallskip
{\it Proof.} 
Since $p\ge1$, the group $\pi_{2l+p}(S^{l+p})$ is either
stable or metastable, so the stable suspension $\Sigma^\infty$ is epimorphic.
Since $l\in\{3,7\}$, the Hopf invariant $H$ is epimorphic.
Hence there is $\varphi\in\pi_{2l+p}(S^{l+p})$ whose stable Hopf 
invariant $H\Sigma^\infty(\varphi)$ is $1\in\Z_2$. 
By the Hopf Invariant Lemma for $q=2l$ and $m=3l+p+1$ we obtain 
$\beta(i\#\overline\varphi)=H\Sigma^\infty(\varphi)=1$. 
\qed

\smallskip
{\bf Realization Lemma.}
{\it Let $N$ be an orientable $(p-1)$-connected $p$-parallelizable closed
$n$-manifold and $n\ge2p+1$. 
Then any homology class $x\in H_p(N;\ \Z\text{ or }\Z_2\ )$ is realizable by an 
embedding 
$\overline x:S^p\times D^{n-p}\to N$.} 

\smallskip
{\it Proof.} Since $N$ is $(p-1)$-connected, any homology class in $H_p(N;\Z)$ 
can be realized by a map $S^p\to N$. 
Hence the same holds for $\Z_2$-coefficients. 
Since $n\ge 2p+1$ every such map is homotopic to an embedding $S^p\to N$. 
Since $N$ is $p$-parallelizable, it follows that this embedding can be extended 
to an embedding $\overline x:S^p\times D^{n-p}\to N$.
\qed

\smallskip
{\bf $\#$-additivity Lemma.}
{\it If $p=3n-2m+3\ge0$, $m\ge n+p+2$ and 
$s:S^p\times D^{n-p}\to N$ is an embedding, then 
$\beta(f\#_sg)=\beta(f)+\beta(g)[s]$ 
for almost embeddings $f:N\to\R^m$ and $g:S^p\times S^{n-p}\to\R^m$,  
where $\beta(g)$ is considered as an element of $\Z_{(m-n-1)}$.} 

\smallskip
{\it Proof.}
Since $m\ge n+p+2$, by the Standardization Lemma (a) we may assume that $f$ 
and $g$ are standardized.
Since $R_m$ and $R_{n-p}$ are isotopic to the identity maps of $\R^m$ and of 
$S^{n-p}$, they do not change orientations. 
Hence 
$$[\Sigma(f\#_sg)]=[\Sigma(f)]+s_*(\id S^p\times R_{n-p})_*[\Sigma(g)],\quad
\text{so we can take}$$ 
$$C_{f\#_sg}:=C_f+s_*(\id S^p\times R_{n-p})_*C_g\quad\text{and}
\quad D_{f\#_sg}:=D_f+R_mD_g.$$ 
We may assume that the supports of $D_f$ and $D_g$ are in $\R^m_+$. 
Identify $S^p\times D^{n-p}$ with $S^p\times D^{n-p}_+$ so that it would 
contain the support of $\widetilde D_g$. 
Then  
$$\beta(f\#_sg)=[\widetilde D_{f\#_sg}]=
[(F\#_sG|_{N-\Int B^n})^{-1}D_f]+
[(F\#_sG|_{N-\Int B^n})^{-1}R_mD_g]=$$
$$=[\widetilde D_f]+s_*[\widetilde D_g]=\beta(f)+[s]\beta(g).\quad\qed$$

%\smallskip
{\it Proof of the New Embedding Theorem.}
Take any $\varphi\in\pi^{m-1}_{eq}(\t N)$.
Since $2m=3n+1-k$, by the Almost Embedding Theorem (c) there is an almost 
embedding $f:N\to\R^m$ such that $\alpha(f)=\varphi$.
By the Realization Lemma there is an embedding $s:S^{k+1}\times D^{n-k-1}\to N$
such that $[s]=-\beta(f)$.
By the Non-triviality Lemma there is an almost embedding 
$G:S^{k+1}\times S^{n-k-1}\to\R^m$ such that $\beta(G)=1$. 
By the $\#$-additivity Lemma $\beta(f\#_sG)=\beta(f)+[s]=0$.
By the completeness of $\beta$-invariant $f\#_sG$ is almost concordant
to an embedding $N\to\R^m$.
(Note that possibly $\alpha(f)\ne\alpha(f\#_sG)$.) 
\qed

%\smallskip

\smallskip
{\it Proof of the New Isotopy Theorem (b).}
Since $2m=3n+2-k$, the surjectivity follows by the Embedding Theorem (b). 
The following proof of the injectivity is a relative version of the proof of 
the New Embedding Theorem.

Take two embeddings $f,f':N\to\R^m$ such that $\alpha(f)=\alpha(f')$.
By the Almost Embeddings Theorem (b) there is an almost concordance $F$ from 
$f$ to $f'$.
By the Realization Lemma there is an embedding $s:S^{k+1}\times D^{n-k-1}\to N$
such that $[s]=-\beta(F)$.

Take an almost embedding $G:S^{k+1}\times S^{n-k}\to\R^{m+1}$ given by the 
Non-triviality Lemma. 
Analogously to the Standardization Lemma we may assume that 

$G(S^{k+1}\times D^{n-k}_-)\subset D^{m+1}_-$ is the standard 
embedding,  

$G(S^{k+1}\times\frac12 D^{n-k}_+)\subset\frac12 D^{m+1}_+$ is the standard 
embedding and 

$G(S^{k+1}\times(D^{n-k}_+-\frac12D^{n-k}_+))
\overset{(*)}\to\subset D^{m+1}_+-\frac12D^{m+1}_+$.  

The latter inclusion (*) gives an almost $G_0$ concordance between standard 
embeddings such that $\beta(G_0)=1$. 
Then $F\#_sG_0$ is an almost concordance from $f$ to $f'$. 

Analogously to the $\#$-additivity Lemma one proves the following.  

%($\#$-additivity for concordances) 
{\it If $p=3n-2m+2\ge0$, $m\ge n+p+3$ and 
$s:S^p\times D^{n-p}\to N$ is an embedding, then 
$\beta(F\#_sG_0)=\beta(F)+\beta(G_0)[s]$ for almost concordances 
$F:N\times I\to\R^m\times I$ and $G_0:S^p\times S^{n-p}\times I\to
\R^m\times I$, 
%N\times I\overset F\to\to\R^m\times I\overset G\to\leftarrow T^n_p\times I,$$ 
where $\beta(G_0)$ is considered as an element of $\Z_{(m-n-1)}$.}   

So $\beta(F\#_sG_0)=\beta(F)+[s]=0$. 
Hence by the completeness of $\beta$-invariant $f$ is isotopic to $f'$.
\qed

%\newpage

\enddocument
\end